\documentclass[11pt,a4paper]{article}
\usepackage[T1]{fontenc}
\usepackage{lmodern,amsmath,amssymb,amsthm,mathtools,microtype}
\usepackage[hmargin=24mm,vmargin=22mm]{geometry}
\usepackage{graphicx}
\usepackage[font=small,labelfont=bf,labelsep=period,hypcap=false]{caption}
\usepackage[hidelinks]{hyperref}
\hypersetup{pdftitle={The Boundary Conjecture for Smoothable Alexandrov Spaces},
 pdfauthor={Nina Lebedeva},
 pdfsubject={Alexandrov geometry, convex hypersurfaces and noncollapsed limits}}
\newtheorem{theorem}{Theorem}
\newtheorem{lemma}[theorem]{Lemma}
\newtheorem{proposition}[theorem]{Proposition}
\theoremstyle{remark}
\newtheorem*{remark}{Remark}
\newcommand{\CBB}{\operatorname{CBB}}
\newcommand{\sn}{\operatorname{sn}_{\kappa}}
\newcommand{\cs}{\operatorname{cs}_{\kappa}}
\newcommand{\md}{m_{\kappa}}
\newcommand{\Int}{\operatorname{Int}}
\newcommand{\dist}{\operatorname{dist}}
\newcommand{\Lip}{\operatorname{Lip}}

\newcommand{\ang}{\angle}
\newcommand{\R}{\mathbb R}
\newcommand{\bd}{\partial}

\begin{document}
\begin{center}
  {\LARGE The Boundary Conjecture for\\[2pt]
  Smoothable Alexandrov Spaces\par}
  \medskip
  {\large Nina Lebedeva\par}
\end{center}
\begin{abstract}
We prove that in a noncollapsed limit of complete Riemannian manifolds
with \(\sec\ge\kappa\), each boundary component of a compact
full-dimensional convex body, equipped with its intrinsic length metric, has
Alexandrov curvature at least \(\kappa\).
We also consider such a body in the limit of a noncollapsed
sequence of Alexandrov spaces with curvature at least \(\kappa\).
If \(\kappa>0\) or the boundary of the body is strongly convex, we show that
the body can be approximated by convex bodies
in the given sequence.
\end{abstract}

\section{Introduction}

A basic conjecture in Alexandrov geometry, stated already by Yuri Burago, Mikhail Gromov, and Grigori Perelman \cite[\S~13.4]{BGP}, asserts that each
component of the boundary of an Alexandrov space, equipped with its intrinsic
length metric, has the same lower curvature bound as the ambient space.

Throughout the paper, we call a subset \emph{convex} if it contains
every minimizing geodesic between any two of its points.
We prove the boundary conjecture for convex bodies in noncollapsed limits of
Riemannian manifolds:
\begin{theorem}\label{main}
Let \(M_j^n\), \(n\ge2\), be complete smooth Riemannian manifolds without
boundary, with \(\sec M_j\ge\kappa\).
Suppose that \(M_j\to A^n\) in the noncollapsed Gromov--Hausdorff sense,
pointed if necessary.
Let \(C\subset A\) be a compact, full-dimensional convex subset
with \(\partial C\ne\varnothing\).

Then each component of \(\partial C\), equipped with its intrinsic length
metric, is \(\CBB(\kappa)\).
\end{theorem}

For convex hypersurfaces in Riemannian manifolds, Sergei Buyalo \cite{Buyalo}
proved the existence of a local lower curvature bound.
Stephanie Alexander, Vitali Kapovitch, and Anton Petrunin \cite{AKP} obtained the optimal
bound.
Karsten Grove, Adam Moreno, and Peter Petersen \cite{GMP} proved the boundary conjecture
for orbit spaces and, more generally, for leaf spaces of singular Riemannian
foliations.

The conjecture has also influenced the development of Alexandrov
geometry. Using the theory of quasigeodesics and gradient curves developed
with Perelman \cite{PP}, Petrunin showed that every shortest path in the
intrinsic boundary metric is a quasigeodesic in the ambient space
\cite[Thm.~1.1]{Pet97}. He later proposed approaching the conjecture through a
Gauss formula for convex hypersurfaces \cite[\S~9]{Pet07}. Nina Lebedeva and Anton Petrunin \cite{LP} constructed a measure-valued curvature tensor for smoothable
Alexandrov spaces as a step in this direction.

In many potential applications, the ambient Alexandrov space arises
as a noncollapsed limit of Riemannian manifolds with a common lower
sectional-curvature bound. Theorem~\ref{main} may already suffice in this setting
(see the graph construction in Section~\ref{sec:further}).

\noindent\emph{Proof strategy.}\enspace
The Riemannian curvature bound applies even to nonsmooth convex
hypersurfaces \cite{AKP}. To prove Theorem~\ref{main}, we therefore lift
suitable convex bodies to the manifolds \(M_j\).
For the full-dimensional bodies considered here, Hausdorff convergence gives
convergence of the corresponding boundary components in their intrinsic
length metrics (Section~\ref{sec:intrinsic}).

Petrunin \cite[Question~9.1.3']{Pet07} asked whether a
convex hypersurface in a noncollapsed Gromov--Hausdorff limit of Alexandrov
spaces can be approximated by convex hypersurfaces in the given spaces.
The following theorem gives the form of lifting that we need.
\begin{theorem}[Lifting of convex surfaces]\label{thm:recovery}\label{recovery}
Let \(X_j^n\) be complete Alexandrov spaces with curvature at least \(\kappa\),
converging without collapse to an \(n\)-dimensional Alexandrov space \(A^n\)
without boundary.
Let \(D\subset A\) be a compact, full-dimensional convex body.
Suppose that either \(\kappa>0\) or \(\partial D\) has a positive lower
strong base-angle bound, as defined in Section~\ref{sec:parallels}.

Then there exist compact, full-dimensional convex bodies
\(D_j\subset X_j\) converging to \(D\) in marked Hausdorff distance along
the prescribed sequence.
\end{theorem}

For \(\kappa\le0\), an arbitrary convex body agrees near each boundary
point with a smaller convex body admitting such approximations;
see Section~\ref{sec:strict}.

In Sections~\ref{sec:parallels}--\ref{sec:cuts}, we construct local
convex supports for inner parallel bodies.
Anton Petrunin \cite[Question~7.3.6]{Pet07} proposed constructing
\((-2+\varepsilon)\)-concave functions, for arbitrarily small
\(\varepsilon>0\), that vanish at a prescribed point and are bounded below
by the negative squared distance from that point. Artem Nepechiy
\cite[Prop.~1.2]{Nep} carried out this construction and obtained functions
that can also be lifted to prescribed noncollapsed approximating sequences.
This construction is what allows us to turn the quadratic margin supplied
by the distance comparison into concave supporting functions.
A suitable linear combination with a modified distance function is concave,
while the quadratic lower bound ensures that its superlevel set contains
the inner parallel body near the point of contact.
A finite family of these supports keeps the convex hulls of interior
approximations close to the corresponding inner parallel body. This gives
the convex approximations in Theorem~\ref{thm:recovery}.

To complete the proof of Theorem~\ref{main}, we apply the lifting theorem
directly to \(C\) when \(\kappa>0\) or the base-angle condition
holds. When \(\kappa\le0\) and \(C\) lacks this bound, we argue locally.
We use strongly concave functions to construct convex bodies \(D_\varepsilon\)
satisfying the required base-angle condition and converging to a body
\(D_0\subset C\) that agrees with \(C\) near a chosen boundary
point \(p\) (Section~\ref{sec:strict}).
We lift each \(D_\varepsilon\) and pass to the limit in the intrinsic
boundary metrics. This gives the curvature bound near \(p\), and
globalization completes the proof (Section~\ref{sec:intrinsic}).

\section{Supporting Functions}\label{sec:parallels}

Throughout Sections~\ref{sec:parallels}--\ref{sec:cuts}, let $X_j\to A$ and $D\subset A$
satisfy the hypotheses of Theorem~\ref{thm:recovery}.

A full-dimensional closed convex subset $D\subset A$, with
the restricted ambient metric, is an Alexandrov space with curvature at least
$\kappa$. Its Alexandrov boundary coincides with its topological boundary in
$A$ \cite[Thm.~1.1]{KZ}, and $D$ is the closure of its interior.

For a body $D$ with $\bd D\ne\varnothing$, write
$f=d_D(\,\cdot\,,\bd D)$ and $D_s=\{f\ge s\}$.
We use the base-angle condition of Stephanie Alexander and Richard Bishop
\cite[\S4]{AB}. A \emph{lower strong base-angle bound}
$\lambda>0$ means that at each $q\in\bd D$,
for every $0<a<\lambda$, all sufficiently short
minimizing chords $qy$, $y\in\bd D\setminus\{q\}$,
satisfy
\begin{equation}\label{ba}
 \ang(\uparrow_q^y,\Sigma_q\bd D)
 \ \ge\ a\,d(q,y)/2 .
\end{equation}
Here $\Sigma_q\bd D$ is the link of the tangent subset
of the boundary in $T_qD$. The radius may depend
on $q,a$. Only boundary points which are nearest
points of interior points will be used in recovery.

\subsection{The quadratic margin}\label{sec:margin}

To construct the desired lifts, we seek local convex supports for the inner
parallel bodies $D_s$ (see Figure~\ref{fig:disks}).
For $p\in\partial D_s$, choose a nearest point $q\in\partial D$.
The open ball $B_s(q)$ is disjoint from $D_s$ and its sphere passes
through $p$, so $d(q,x)\ge d(q,p)=s$ for every $x\in D_s$.

The idea is to combine a modified distance function
with a strongly concave correction. The distance term must separate $D_s$
from the point of contact by a positive quadratic amount, while the correction
must be strong enough to make the resulting function concave. The purpose of
the following comparison is to show that these two requirements are compatible.

Write \(\sn(r)=\sin(\sqrt\kappa\,r)/\sqrt\kappa\) for \(\kappa>0\),
\(\sn(r)=r\) for \(\kappa=0\), and
\(\sn(r)=\sinh(\sqrt{-\kappa}\,r)/\sqrt{-\kappa}\) for \(\kappa<0\).
Set \(\cs=\sn'\) and \(\md(r)=\int_0^r\sn(v)\,dv\).
Write \(d_q=d(q,\cdot)\).
Along unit-speed minimizing segments the following comparison
inequality holds across the cut locus \cite[(1.1)]{AB},
\begin{equation}\label{upper}
  (\md\circ d_q)''\le \cs(d_q).
\end{equation}
For \(\kappa>0\) we take \(0<s<\pi/(4\sqrt\kappa)\).

\begin{lemma}[Two distances and one foot]\label{foot}
Fix \(p\) with \(f(p)=s>0\), and fix \emph{any} \(q\in\bd D\) satisfying
\(d(p,q)=s\). Suppose \eqref{ba} holds at \(q\) with coefficient \(a\).
The choice \(a=0\) requires no extra hypothesis. Then
\begin{equation}\label{two}
 \liminf_{\substack{x\to p,\ x\in D_s\\x\ne p}}
 \frac{\md(d(q,x))-\md(s)}{d(p,x)^2}
 \ \ge\ L_s,\qquad
 L_s=1-\tfrac12\cs(s)+\tfrac a2\sn(s).
\end{equation}
\end{lemma}

\begin{proof}
Let \(x\in D_s\setminus\{p\}\) tend to \(p\). Write \(r=d(p,x)\) and \(\ell=d(q,x)\). The perimeter of \(pqx\) is
at most \(2s+2r\). Since \(x\in D_s\),
\(\ell\ge s\). It suffices to treat a sequence on which the quotient
on the left of \eqref{two} is bounded above. Then \(\ell-s=O(r^2)\)
as \(\md'(s)>0\). Let \(\beta\) be the comparison angle at \(q\) in
\(pqx\). The model cosine law reads
\[
 \md(r)=\md(\ell-s)+\sn(s)\sn(\ell)(1-\cos\beta).
\]
It follows that \(\beta\to0\) and \(\sn(s)\sin\beta/r\to1\).
In particular \(\beta>0\) eventually.

Every path between the two copies of \(p\) in the double meets
\(\bd D\), so the path through the chosen foot \(q\) minimizes and
has length \(2s\). This gives \(T_qD=\R_+\times T_q\bd D\), with
\(qp\) along the normal factor \cite[proof of Thm.~3.3.1, note 19]{Pet07}.
For any chosen \(qx\), the splitting gives \(B=\ang pqx\le\pi/2\).
Hinge comparison gives \(B\ge\beta\).
Pass to a subsequence so that the normalized tangential projections of
these directions tend to a unit vector \(\xi\). To realize \(\xi\)
by points of \(\bd D\), take a unit-speed quasigeodesic in this
extremal subset starting at \(q\) in direction \(\xi\)
\cite[\S5.1, properties 1, 2 and 4]{Pet07}. At parameter \(u\) its
distance to \(q\) is \(u+o(u)\). By continuity it contains a point \(y\)
at distance exactly \(r\), with parameter between \(r/2\) and \(2r\).
The directions \(qy\) tend to \(\xi\). Their elevations \(\phi\) above
the tangential factor tend to zero, and the angle \(\delta\) between
the two tangential projections tends to zero. By \eqref{ba},
\(\phi\ge ar/2\).

For \(\theta=\ang xqy\) the product metric gives
\[
 \cos\theta=\cos B\sin\phi+\sin B\cos\phi\cos\delta.
\]
As a function of \(B\) on \([\beta,\pi/2]\) the right hand side is
concave for all sufficiently small \(r\). Its values at the endpoints
tend to \(0\) and \(1\), respectively, so its minimum occurs at
\(\beta\). At \(B=\beta\) it is increasing in \(\phi\) between \(ar/2\)
and the actual elevation, since \(\beta,\phi\to0\). Substitution yields
\begin{equation}\label{angle}
 \cos\theta\ge \sin\beta+ar/2+o(r).
\end{equation}
In fact the angular mismatch is multiplied by \(\sin\beta=O(r)\).
Thus \(\delta\to0\) suffices, and no bound on the rate at which \(B\)
approaches a limit is involved.

Now \(s\le f(x)\le d(x,y)\). Since \(d(q,y)=r\), the perimeter of
\(qxy\) is at most \(2s+4r\). For \(\kappa>0\) both perimeters are
therefore below \(2\pi/\sqrt\kappa\) for small \(r\), and the side
lengths lie in the increasing range of \(\md\).
Hinge comparison for \(qx,qy\) and the
model cosine law give
\[
 \md(s)\le\md(d(x,y))
 \le \md(\ell)+\cs(\ell)\md(r)-\sn(\ell)\sn(r)\cos\theta .
\]
Insert \eqref{angle}, use \(\md(r)=r^2/2+O(r^4)\) and
\(\sn(r)=r+O(r^3)\), and divide by \(r^2\). The resulting lower
bound is \(L_s\). Choosing a subsequence attaining any putative
smaller liminf justifies the directional subsequence used above.
\end{proof}

If the strong base-angle hypothesis holds with bound \(\lambda>0\),
fix \(0<a<\lambda\). Otherwise \(\kappa>0\) and we take \(a=0\). There is \(s_*>0\), smaller than
\(\max f\), such that for every \(0<s<s_*\)
\begin{equation}\label{gap}
 2L_s-\cs(s)=2(1-\cs(s))+a\sn(s)>0,\qquad \cs(s)>0.
\end{equation}
For positive curvature \(1-\cos(\sqrt\kappa\,s)>0\). For \(\kappa=0\)
the gap is \(as\). For \(\kappa=-k^2<0\) it is positive exactly when
\(a>2k\tanh(ks/2)\). In particular, only shallow levels are needed.

\begin{lemma}\label{shallow}
For \(0<t<s_*\), the inner parallel body \(D_t\) is convex.
\end{lemma}
\begin{proof}
By Petrunin's theorem on the isometry of tangent cones
along geodesics \cite[Thm.~1.1(A)]{Pet98},
a minimizing geodesic in $D$ meeting $\bd D$
at an interior point lies entirely in $\bd D$.
In particular, a chord with one interior endpoint
has its open segment in $\Int D$.

A minimizing chord between points of \(D_t\) stays in \(D\) by convexity.
Since both endpoints have depth at least \(t>0\), the whole chord has
positive depth. If its
minimum depth were \(s<t\), it would be attained at an interior point
\(p=\gamma(0)\) of a unit-speed parametrization. Fix a foot \(q\) of
\(p\). Both points \(\gamma(\pm h)\) belong to \(D_s\).
By \eqref{two} the lower limit of the centered second difference of
\(\md(d_q\circ\gamma)\), divided by \(h^2\), is at least \(2L_s\).
By \eqref{upper} its upper limit is at most \(\cs(s)\).
This contradicts \eqref{gap}.
\end{proof}

\subsection{Construction of concave supports}\label{sec:supports}

\begin{center}
\refstepcounter{figure}\label{fig:disks}
\includegraphics[width=.94\linewidth]{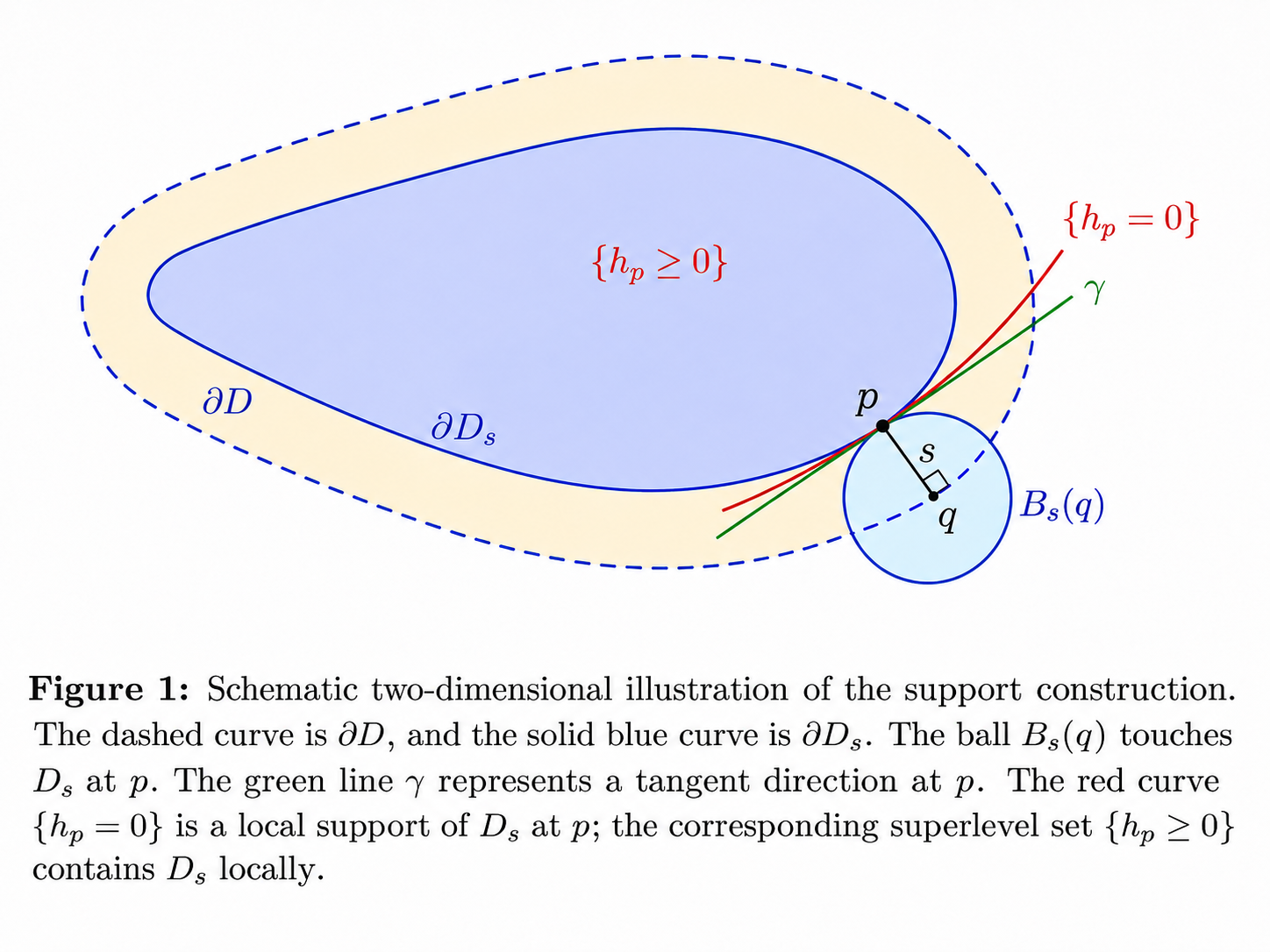}
\end{center}

We now turn to the construction of supporting functions, using the quadratic
margin obtained in the preceding subsection.

There are two difficulties in constructing such supports: their superlevel
sets must be convex, and they must touch the body without cutting into it.

Distance functions are natural building blocks for semiconcave functions in
Alexandrov spaces, but their superlevel sets need not be convex. Moreover, a
ball touching a convex body may still cut into it nearby, so the corresponding
distance superlevel need not support the body.

This is why we construct supports for the inner parallel body $D_s$.
By definition, $D_s$ is disjoint from every open ball of radius $s$
centered on $\partial D$, and each point of $\partial D_s$ lies on at
least one of the corresponding spheres. Thus at each such point we
have a distance superlevel containing $D_s$ and touching it there.
It remains to make the superlevel convex while preserving this support.

Perelman's construction
\cite[\S3.6]{PerMorse} provides strongly concave functions whose first-order
behavior at a chosen point can be prescribed up to an arbitrarily small error.
Their superlevel sets, however, need not support a given body either: they may
bend too sharply and cut into it. Answering a question of Petrunin
\cite[Question~7.3.6]{Pet07}, Nepechiy \cite[Prop.~1.2]{Nep} constructed
\((-2+\eta)\)-concave functions, for arbitrarily small \(\eta>0\), bounded
below by the negative squared distance.
Unlike the first-order control provided by Perelman's construction, this gives
a quantitative bound on how far the function can fall below its maximum.
This quadratic control is precisely what we need to preserve the margin
obtained in Section~\ref{sec:margin}. Nepechiy's functions, however, have a maximum at their
center and do not themselves determine a supporting direction.

Roughly speaking, the distance term provides the contact and the inward
direction of the support, while Nepechiy's function corrects the lack of
convexity of the distance superlevel. The quadratic margin established in
Section~\ref{sec:margin} ensures that this correction does not cut into $D_s$.

\noindent\textbf{Nepechiy's construction.}
Nepechiy's Theorem B \cite{Nep} provides \((-2+\eta)\)-concave functions
with a quadratic lower bound, for arbitrarily small \(\eta>0\).
Proposition~1.2 (case $k=0$) strengthens this result by
producing functions that are $\eta$-constructible in the sense of
\cite[Definition~1.1]{Nep}. In particular, these functions can be lifted to
prescribed noncollapsed approximating sequences while retaining the estimates
needed below.

Let \(A\) be an \(n\)-dimensional Alexandrov space without boundary,
and let \(p\in A\). For every \(\eta>0\), there exist \(R>0\) and
a function \(b_p\) on \(B_R(p)\) such that
\[
b_p(p)=0,\qquad \operatorname{Lip}b_p\le 2R,\qquad
b_p''\le-2+\eta,\qquad
-d(p,x)^2\le b_p(x)\le0.
\]

Moreover, for any prescribed noncollapsed sequence of
\(n\)-dimensional Alexandrov spaces with a common lower curvature bound,
converging to \(A\), and any choice of points \(p_j\to p\), there exist,
for all sufficiently large \(j\), functions \(b_{p,j}\) on \(B_R(p_j)\)
satisfying
\[
b_{p,j}(p_j)=0,\qquad
\operatorname{Lip}b_{p,j}\le2R,\qquad
b_{p,j}''\le-2+\eta,
\]
and
\[
b_{p,j}(x)\ge-d(p_j,x)^2
\quad\text{for }d(p_j,x)\ge\varepsilon_jR,
\qquad \varepsilon_j\to0.
\]
Fix \(0<s<s_*\) and a point \(p\) with \(f(p)=s\).
Let \(q\in\bd D\) be any nearest boundary point to \(p\).
By \eqref{gap}, choose
\[
 \cs(s)<\alpha<2L_s,\qquad
 0<\eta<1,\qquad \alpha(1-\eta/2)>\cs(s),\qquad
 \alpha/2<d_0<L_s .
\]
For this choice of \(\eta\), take \(b_p\) as above and set
\begin{equation}\label{support}
  h_p=m_\kappa(d_q)-m_\kappa(s)+\tfrac{\alpha}{2}b_p .
\end{equation}
By \eqref{upper} and the $(-2+\eta)$-concavity of $b_p$, on a sufficiently small ball about $p$,
\[
h_p''\le \cs(d_q)-\alpha(1-\eta/2)<0.
\]
On the other hand, \eqref{two} and $b_p\ge-d_p^2$ give, after shrinking the ball,
\[
h_p(x)\ge (d_0-\alpha/2)d(p,x)^2>0
\qquad (x\in D_s\setminus\{p\}).
\]

Here \(h_p(p)=0\), and \(h_p<0\) on \(pq\) immediately toward \(q\),
where \(d_q<s\). Every minimizer with both endpoints sufficiently
close to \(p\) stays in the ball of concavity by its length. Thus the
superlevel is locally convex and supports \(D_s\) at \(p\).

\begin{remark}[Euclidean interpretation]
In Euclidean space, the difference of two squared distance functions
is affine. In an Alexandrov space, the analogous difference is locally
a difference of semiconcave functions, but need not
be either convex or concave.
Our construction can be viewed as replacing the negative
squared-distance term by Nepechiy's almost-quadratic concave function.
\end{remark}

\section{Finite cuts on the original sequence}\label{sec:cuts}

The closed convex hull of a subset $E\subset X$
is the intersection of all closed convex
subsets of $X$ containing $E$.

Fix a depth $0<t<s_*$. Consider sets $E_j\subset X_j$ whose images
under the Gromov--Hausdorff approximations $F_j$ lie within $o(1)$
of $D_t$. To control their hulls in $X_j$, we lift finitely many local
supports from the boundary of the larger body $D_s$, where $0<s<t$.
For all large $j$, the cuts give a compact convex
set $K_j\subset X_j$ with
\[
 E_j\subset K_j,\qquad F_j(K_j)\subset D_s.
\]
The entire hull of $E_j$ lies in $K_j$ by definition. Its marked image
therefore lies in $D_s$, and the compact sets $D_s$ decrease to $D_t$
as $s\uparrow t$.

\begin{lemma}[Convex hulls at a fixed depth]\label{lem:hull-control}
Suppose $\bd D\ne\varnothing$ and $0<t<s_*$.
Choose marked base points $o_j\to o$, a radius $R$ with
$D\subset B_{R/2}(o)$, and Gromov--Hausdorff marks $F_j$ on
$\overline B_{4R}(o_j)$ with distortion $\epsilon_j\to0$.
Let $E_j\subset B_R(o_j)$ be compact sets with
\[
 \sup_{e\in E_j}\dist_A(F_je,D_t)\longrightarrow0.
\]
Then their closed convex hulls $H_j$ are eventually compact,
contained in $B_R(o_j)$, and satisfy
\[
 \sup_{x\in H_j}\dist_A(F_jx,D_t)\longrightarrow0.
\]
If $E_j$ converge to $D_t$ in marked Hausdorff distance, then so do $H_j$.
\end{lemma}

\begin{proof}
Fix $0<s<t$.

\smallskip\noindent\textbf{Lifted supports.}
Extend $f$ to $A$ by $f(x)=d_A(x,A\setminus D)$.
This is 1-Lipschitz, agrees with the depth on $D$, and vanishes outside.
On the marked ball put
\[
 u_j(x)=\inf_z\{f(F_jz)+d_j(x,z)\},\qquad
 f(F_jx)-\epsilon_j\le u_j(x)\le f(F_jx).
\]
The infimum is over the marked ball. It gives a continuous function even
for discontinuous marks. Its $s$-superlevel there lies in $B_R(o_j)$
for all large $j$. A minimizer with both endpoints in that smaller ball
is contained in $B_{3R}(o_j)$.

For this $s$, choose $\alpha,\eta,d_0$ as in
Section~\ref{sec:supports}. For each $p$ with $f(p)=s$,
choose one foot $q$ and lifts $p_j,q_j$.
Applying the constructibility statement recalled in
Section~\ref{sec:supports}, we obtain functions $b_j$
on balls $B_{R_p}(p_j)$, where $R_p>0$ is independent
of $j$, satisfying the estimates stated there.
In particular, the quadratic lower bound holds
outside $B_{e_{p,j}R_p}(p_j)$, where $e_{p,j}\to0$.
Define
\[
 h_{p,j}=\md(d_j(q_j,\cdot))-\md(s)+(\alpha/2)b_j.
\]
We need these functions to be small near $p_j$ and uniformly positive
on a fixed annulus where $u_j\ge s$. A positive cutting level can then
remove a neighborhood of $p_j$ without reaching the radial boundary.
The normalization gives $h_{p,j}(p_j)\to0$, and the Lipschitz constants
are bounded on a fixed ball. By \eqref{upper},
$h_{p,j}''\le\cs(d_j(q_j,\cdot))-\alpha(1-\eta/2)<0$
on $B_{2r_p}(p_j)$ after choosing $2r_p<R_p$ sufficiently small,
then $j$ sufficiently large.

For the annular lower bound, use \eqref{two} to decrease $r_p$ until
the numerator there is at least $d_0d(p,x)^2$ for
$x\in D_s\cap B_{2r_p}(p)$. Also require $r_p<(t-s)/4$ and $2r_p<s$.
On $r_p/2\le d_j(x,p_j)\le r_p$ the excluded ball is eventually
absent. At the points where $u_j(x)\ge s$, the inclusion
$F_jx\in D_s$ gives, uniformly on this fixed annulus,
\[
 h_{p,j}(x)\ge(d_0-\alpha/2)d(p,F_jx)^2-o_j(1)\ge4a_p
\]
for a fixed $a_p>0$ and all large $j$. Indeed
$d(p,F_jx)\ge r_p/3$ eventually.

Choose $0<\tau_p<r_p/10$ so that $h_{p,j}<a_p/2$ on
$B_{\tau_p}(p_j)$, using the Lipschitz bound and the center value.
Take a finite subcover of $\{f=s\}$ by $B_{\tau_p/3}(p)$, with centers
$p_i$, before making the final choice of $j$.

\smallskip\noindent\textbf{The finite cuts.}
On the marked ball form the closed set
\[
 K_j=\{u_j\ge s\}\setminus\bigcup_i V_{i,j},\qquad
 V_{i,j}=\{d_j(x,p_{i,j})<r_i,\ h_{i,j}(x)<a_i\}.
\]
Every point of $E_j$ has marked depth at least $t-o(1)$ and distance
at least $t-s-o(1)$ from each center, uniformly over $E_j$.
Thus $E_j\subset K_j$.
The sets $K_j$ are compact, being closed inside the marked compact
ball and contained in its interior.

The finite cover also gives a fixed $\delta>0$, with $s+\delta<s_*$,
such that $u_j\ge s+\delta$ on $K_j$ eventually.
Otherwise points of $K_j$ of depth tending to $s$ have a marked limit
on $\{f=s\}$, in a member of the finite cover. They eventually lie in
the corresponding removed open set.

\smallskip\noindent\textbf{Convexity.}
We first show that every minimizing chord with endpoints in $K_j$
stays in $\{u_j>s\}$.
If such chords could reach $u_j\le s$ for infinitely many $j$,
choose one offending chord for each such $j$. Their lengths are bounded.
After passing to a subsequence they converge to a minimizer with endpoints
in $D_{s+\delta}$ meeting $\{f\le s\}$, contrary to Lemma~\ref{shallow}.

Such a chord $\gamma$ cannot enter any $V_{i,j}$ either.
If it does, take a connected component of $\gamma^{-1}(V_{i,j})$.
At a component endpoint the chord cannot meet the radial wall
$d_j(x,p_{i,j})=r_i$, since $h_{i,j}\ge4a_i$ there.
Both endpoints therefore have $h_{i,j}=a_i$, including when an endpoint
is a chord endpoint. The intervening subsegment lies in the domain of
concavity and has $h_{i,j}<a_i$, a contradiction.
Thus $K_j$ is convex in $X_j$.

By the intersection definition of the hull,
\[
 E_j\subset H_j\subset K_j,\qquad F_j(K_j)\subset D_s.
\]
In particular, $H_j$ is compact. For each fixed $s<t$ these inclusions
hold eventually. The decreasing compact sets $D_s$ converge to $D_t$
as $s\uparrow t$, which gives the asserted estimate.
The reverse Hausdorff inclusion follows from the sets $E_j$ when they
converge to $D_t$.
\end{proof}

\begin{proof}[Proof of Theorem~\ref{thm:recovery}]
If \(\bd D=\varnothing\), then \(D=A\) is compact. Choose
\(R>\operatorname{diam}A+1\). Pointed convergence rules out points
at distance \(R\) from the marked basepoints in \(X_j\) for all
large \(j\); a minimizing segment to a more distant point would
supply such a point. Hence \(X_j\) is compact and converges globally
to \(A\), and we take \(D_j=X_j\) for all sufficiently large \(j\).
The finitely many initial indices are handled by the superlevel bodies
constructed below. We therefore assume \(\bd D\ne\varnothing\).

Take \(t_m\downarrow0\), with $t_m<s_*$.
For fixed \(m\) lift increasingly fine nets of \(D_{t_m}\) and add a small
closed source ball centered at a lift of an interior point of depth
greater than \(t_m\). Its radius can tend to zero.
By Lemma~\ref{lem:hull-control}, their closed convex hulls
are eventually compact and converge to \(D_{t_m}\). The added balls
make the hulls full-dimensional.
By density of \(\Int_A D\), \(D_{t_m}\to D\).

Choose increasing indices \(J_m\) so that the hull error at level
\(t_m\) is \(<1/m\) for \emph{every} \(j\ge J_m\). Use level \(t_m\)
on \(J_m\le j<J_{m+1}\). This gives recovery along the whole tail.
A lift of a fixed point outside the closed body \(D\) lies outside
these hulls for large \(j\), by marked Hausdorff convergence.

At each of the finitely many earlier indices choose
\(z\in\Int X_j\) and a strongly concave function \(b\) on \(B_r(z)\)
with an isolated maximum \(b(z)=0\), by
\cite[Thm.~B and Cor.~2.3]{Nep}. Shrink \(r\) so that \(b<0\)
off \(z\). Choose \(c>0\) with \(b<-c\) on
\(r/4\le d_j(z,\cdot)\le r/2\). Then
\[
 Q=\{x\in\overline B_{r/2}(z):b(x)\ge-c\}\subset B_{r/4}(z)
\]
is compact and full-dimensional. Every minimizing chord between its
points stays in \(B_r(z)\) by its length; concavity and the annular
inequality keep it in \(Q\). Use these bodies at the earlier indices.
\end{proof}

\section{Making a local boundary piece strictly convex}\label{sec:strict}

In Theorem~\ref{main}, the limit \(A\) has empty boundary by \cite[Thm.~6.4]{DK}, since the
complete source manifolds have empty boundary and the convergence
preserves dimension.

\begin{proposition}\label{strict}
For each \(p\in\bd C\) there is a compact full-dimensional convex body
\(D_0\subset C\) that agrees with \(C\) near \(p\).
The body \(D_0\) is the increasing Hausdorff limit of compact convex
bodies \(D_\varepsilon\) with positive lower strong base-angle bounds.
\end{proposition}

\begin{proof}
Take a Lipschitz function \(\psi\) with \(\psi(p)=0\), an isolated
maximum at \(p\), and \(\psi''\le-b<0\) on \(B_{2R}(p)\)
\cite[Prop.~1.2(i)--(ii)]{Nep}. Shrink \(R\) so \(\psi<0\) off \(p\).
On the compact annulus \(R/8\le d(p,\cdot)\le R/2\) its maximum is
negative. Choose \(c>0\) small so that
\[
 Q=\{x\in\overline B_{R/2}(p):\psi(x)\ge-c\}\ \subset B_{R/8}(p).
\]
Every minimizing chord with endpoints in \(Q\) stays in \(B_{R/2}(p)\)
by its length, and stays in \(Q\) by concavity. In particular \(Q\) is
compact, convex and contains an ambient neighborhood of \(p\).

Put \(K=\max(0,-\kappa)\). Reduce \(c\) if needed and choose
\(\tau>0\) with \(b>K(\tau+c)\). Write \(f=d_C(\cdot,\bd C)\).
On \(\Int C\) one has \(f''\le Kf\) along unit-speed geodesics.
For \(K=0\) this is boundary-distance concavity. For \(K>0\),
\cite[Thm.~3.3.1]{Pet07} gives
\((\sinh(\sqrt K f))''\le K\sinh(\sqrt K f)\).
Writing \(k=\sqrt K\), the one-dimensional chain rule gives
\(D^2(\sinh(kf))=k\cosh(kf)D^2f+k^2\sinh(kf)(f')^2\,dt\).
Division by the positive coefficient of \(D^2 f\) gives
\(f''\le\sqrt K\tanh(\sqrt K f)(1-(f')^2)\le Kf\).
These inequalities hold weakly along each geodesic.

On \(\Int C\cap B_R(p)\) define
\[
 H_\varepsilon=\varepsilon\psi+\min\{f-\varepsilon\tau,\varepsilon c\},
 \qquad D_\varepsilon
 =\{x\in C\cap Q:f(x)\ge\varepsilon(\tau-\psi(x))\}.
\]
In \(B_{R/2}(p)\cap\Int C\) this is the superlevel \(H_\varepsilon\ge0\).
On \(D_\varepsilon\) we have \(f\ge\varepsilon\tau>0\).
The first branch can be active only where \(f\le\varepsilon(\tau+c)\).
There its second derivative after adding \(\varepsilon\psi\) is at most
\(-\varepsilon[b-K(\tau+c)]\). The other branch has bound
\(-\varepsilon b\). At equality the first bound persists locally with
half its margin, by continuity of \(f\). Where the first branch is
strictly inactive only the second matters. A minimum of functions with
the same concavity bound retains it. Thus locally, and hence along
each chord in the working domain,
\[
 H_\varepsilon''\le-b_\varepsilon,\qquad
 b_\varepsilon=\tfrac{\varepsilon}{2}[b-K(\tau+c)]>0 .
\]
All minimizing chords of \(D_\varepsilon\) stay in \(Q\) and in
\(\Int C\), so they stay in this superlevel. Fix \(z\in\Int C\) near
\(p\) with \(\psi(z)>-c\). A small fixed ball about \(z\) satisfies
both defining inequalities strictly for all sufficiently small
\(\varepsilon\). Thus the bodies are full-dimensional.
Every actual boundary point has \(H_\varepsilon=0\), since the outer
ball constraint is inactive. Conversely a zero in the body admits a
direction of positive differential towards \(z\). On the open ambient
domain a supporting vector \(v\) exists with \(dH(w)\le-\langle v,w\rangle\)
for all \(w\) \cite[Lemma~1.3.7]{Pet07}.
It is nonzero, and \(dH(v)<0\). A zero therefore cannot lie in the
interior of the superlevel.

For a boundary chord of length \(\ell\) with initial unit vector \(u\)
at \(x\), strong concavity on the chord implies
\(d_xH_\varepsilon(u)\ge b_\varepsilon\ell/2\).
For every \(\zeta\in\Sigma_x\bd D_\varepsilon\) the same differential
vanishes at \(\zeta\), since boundary secants converge to it in a blow-up of
the boundary pair \cite[Cor.~6.5]{DK}, and \(H_\varepsilon\) vanishes
at all those secants' endpoints. Their rescaled difference quotients
are therefore zero \cite[Def.~1.3.1]{Pet07}. The isometric inclusion
\(D_\varepsilon\subset A\) identifies its tangent cone with a subcone
of \(T_xA\). The blow-up differential is \(L_\varepsilon\)-Lipschitz,
\(L_\varepsilon=1+\varepsilon\Lip\psi\).
For unit vectors the cone distance is \(2\sin(\ang(u,\zeta)/2)\le
\ang(u,\zeta)\). Taking the infimum over \(\zeta\) proves \eqref{ba}
with \(a=b_\varepsilon/L_\varepsilon>0\), uniformly in \(x\).

Set \(D_0=C\cap Q\). The \(D_\varepsilon\) increase as
\(\varepsilon\downarrow0\). From any \(x\in D_0\), the open chord to
\(z\) is in \(\Int C\) and has \(\psi>-c\). Each of its points lies in
\(D_\varepsilon\) for all sufficiently small \(\varepsilon\).
This proves density and hence Hausdorff convergence to \(D_0\).
The previously chosen inball works throughout. Since \(\psi(p)=0\),
the cut defining \(Q\) is inactive near \(p\).
\end{proof}

\begin{remark}
Even smooth marks close to the identity in \(C^1\) can destroy the
convexity of a graph. On the Euclidean plane consider
\(F_j(u,v)=(u,v+j^{-2}\sin(ju))\). The image of the epigraph of
\(bu^2/2\), \(0<b<1\), has lower graph with second derivative
\(b-\sin(ju)\), which is negative at \(u=\pi/(2j)\).
The cuts in Theorem~\ref{thm:recovery} obtain convexity directly in the
source metric.
\end{remark}

\section{Proof of Theorem~\ref{main}}\label{sec:intrinsic}

Under noncollapsed Gromov--Hausdorff convergence of compact Alexandrov
spaces with a common lower curvature bound, corresponding boundary
components converge in their intrinsic length metrics, provided the
limit has nonempty boundary. Indeed, the boundaries converge as
extremal subsets \cite[Thms.~6.4 and 6.8]{DK}, their components
correspond by relative stability \cite[Thm.~9.2]{Kap}, and their
intrinsic metrics converge by \cite[Thm.~1.2]{Pet97}.
We use this for Hausdorff-convergent full-dimensional convex bodies
with the restricted ambient metrics; their Alexandrov boundaries
coincide with their topological boundaries in the surrounding
boundaryless spaces \cite[Thm.~1.1]{KZ}.

For a body recovered by Theorem~\ref{thm:recovery}, applied with \(X_j=M_j\), each approximating
boundary is locally a convex hypersurface of the smooth \(M_j\).
Indeed, in a strongly convex normal neighborhood of a boundary point
\(x\), take rays from a nearby interior point through a neighborhood
of \(x\). On the ray through \(x\), points just before \(x\) are
interior and points just beyond are exterior by the open-chord
property. These two strict conditions persist on nearby rays.
Convexity gives a unique exit on each, and the exit time is continuous
by openness of the interior and exterior. This gives the local graph.
The theorem of Alexander--Kapovitch--Petrunin \cite[Thm.~1]{AKP}
gives intrinsic curvature at least \(\kappa\), including for nonsmooth
convex hypersurfaces. By intrinsic convergence of the boundary components and
GH closedness of four-point comparison, the recovered body's boundary
has the same bound. This already proves the theorem when \(\kappa>0\).

For \(\kappa\le0\) apply this argument first to each fixed
\(D_\varepsilon\) of Proposition~\ref{strict}. Then send
\(\varepsilon\downarrow0\), using intrinsic convergence of boundary components
for \(D_\varepsilon\to D_0\). Hence each intrinsic boundary component of \(D_0\) is \(\CBB(\kappa)\). The two boundary sets \(\bd D_0\) and
\(\bd C\) agree in some \(B_\rho(p)\). Their intrinsic metric germs
agree as well. Indeed, the intrinsic topology on either extremal
boundary is its subspace topology \cite[Lemma~9.8]{Kap}. In a sufficiently
small common neighborhood of \(p\), the intrinsic distances to \(p\)
in both boundaries are \(<\rho/8\) and the ambient distances are
\(<\rho/8\). Between two such points paths of length \(<\rho/2\)
are available in either boundary, whereas leaving \(B_\rho(p)\)
costs at least \(7\rho/8\). Taking infima over paths in the common
set therefore gives the same distances.

Each component of \(\bd C\) is compact in its intrinsic topology, hence
a complete length space. It is locally \(\CBB(\kappa)\) at every point.
The globalization theorem \cite[Thm.~8.31]{Found} proves the assertion.

When $n=2$ and $\kappa>0$, each boundary component has length at most
$2\pi/\sqrt{\kappa}$ by \cite[\S3.3.5]{Pet07}, so the conclusion also
holds under the convention imposing the diameter bound in dimension one.
\par\hfill$\square$\par

\section{Further remarks}\label{sec:further}

\smallskip\noindent\textbf{Lifting and varying functions.}
Lifting a function is a natural counterpart to lifting a convex body:
the epigraph of a locally convex function, or the hypograph of a locally
concave one, is locally convex in $A\times\mathbb R$.
After local strictification, one can apply the lifting construction to
compact portions of these sets. Reading off their vertical fibers on
smaller domains gives functions converging to the original one.
Varying the epi- or hypographs themselves, while keeping their vertical
sections as rays, similarly gives variations of the functions.
For a $\lambda$-concave function one uses a parabolic warped product,
as in \cite[p.~159]{Pet07}, with a change of height coordinate preserving
$\lambda$. This gives a local answer to \cite[Question~9.1.3]{Pet07}.
On a compact boundaryless limit, the local liftings can also be patched
to give global liftings: for $\lambda>0$, a finite minimum gives an
arbitrarily small loss of concavity, which is removed by rescaling.
The case $\lambda=0$ is constant.
The lifted functions also approximate the local gradient flow of the
original function, by \cite[\S2.2]{Pet07}.

\smallskip\noindent\textbf{Variations of convex bodies.}
The construction with finite supporting cuts in Section~3 also allows
one to vary convex bodies in Alexandrov spaces, without assuming
smoothability. By choosing different support points on the inner parallel
bodies and independently varying the cutting levels, one obtains families
of convex bodies between $D_t$ and $D_s$, for small $0<s<t$. The cuts
can be made in disjoint neighborhoods, so that different parts of the
boundary can be varied independently.

\smallskip\noindent\textbf{Deforming limits and their approximations.}
Alexandrov spaces also arise in Riemannian geometry as limits of
sequences of possible counterexamples to conjectured statements.
Varying such a limit may help us understand nearby geometries.
Moreover, in the smoothable case, the construction described below
allows us to vary not only the limit but also a prescribed approximating
sequence of manifolds.
Concerning the boundary conjecture, Petrunin suggested varying
Alexandrov spaces in the following way \cite{Pet07}.
For a concave function $f$ on a convex neighborhood $U\subset A$,
consider the graphs
\[
 \Gamma_t=\{(x,tf(x)):x\in U\}\subset A\times\mathbb R,
 \qquad t\ge0,
\]
with their intrinsic metrics.
If $A$ is smoothable, Theorem~1 gives these graphs local
curvature at least $\kappa_0=\min\{\kappa,0\}$, and they converge
locally to $U$ as $t\to0$.
This construction can also be carried out along a prescribed smoothing
sequence $M_j\to A$: lift $f$, add a vanishing strictly concave
correction, and smooth on smaller domains using \cite[Lemma~4]{AKP}.
For each fixed $t$, the resulting smooth graph metrics converge to the
intrinsic metric of $\Gamma_t$, with the same lower curvature bound
$\kappa_0$.
On compact spaces one can instead use semiconcave functions and the
parabolic construction. Provided these functions admit lifts along the
prescribed sequence with uniform semiconcavity and Lipschitz bounds,
the resulting graph metrics have curvature at least
$\kappa-O(t^2)$. For $\kappa\ne0$, a homothety tending to the
identity then restores the lower bound $\kappa$.

\smallskip\noindent\textbf{Stability of convex hulls.}
Our approach could also be viewed in terms of the stability of convex
hulls under Gromov--Hausdorff convergence. One lifts finite nets
from the limiting convex body and takes their convex hulls in the
approximating spaces. The main difficulty is to keep the whole hull close
to the body one wants to recover. Alexander Lytchak and Anton Petrunin \cite{LytPetConvex} show that, for a
\(C^2\)-generic Riemannian metric in dimension at least three, connected
convex sets have interior unless they are contained in a geodesic.
Under the hypotheses of Theorem~\ref{thm:recovery}, our supporting cuts
provide the required control for inner approximations;
see Lemma~\ref{lem:hull-control}.

If the additional hypothesis in Theorem~\ref{thm:recovery} is dropped,
such stability may fail even at \(\kappa=0\). For example, a narrow central
belt on the intrinsic boundary surface of a long solid cylinder is convex,
but admits no compact convex recovery on slightly tapering capped frusta.

\section*{Acknowledgments and disclosure}

I am grateful to Anton Petrunin for introducing me to Alexandrov geometry and, in particular, to the boundary conjecture, and for many discussions over the years.

This paper was developed and written in collaboration with ChatGPT (OpenAI).
I fully understand all mathematical arguments presented here and take full responsibility for its content.

\end{document}